# Harnessing GPU Acceleration in Large-Scale Process Optimization

**Boxun Huang[a], David Y. Shu[ab], Michel Schanen[c], Mihai Anitescu[c], Rahul Gandhi[d], and Sungho Shin[a*]**

[a] Massachusetts Institute of Technology, Department of Chemical Engineering, Cambridge, MA 02139, USA
[b] Massachusetts Institute of Technology, MIT Energy Initiative, Cambridge, MA 02139, USA
[c] Mathematics and Computer Science Division, Argonne National Laboratory, Lemont, IL 60439, USA
[d] Shell Technology Center, Houston, TX 77082, USA
[*] Corresponding Author: sushin@mit.edu

## ABSTRACT

This paper presents a proof-of-concept workflow for equation-oriented process optimization that runs entirely on a GPU. Process optimization models often incorporate complex interconnected unit operations, dynamics, and uncertainties, resulting in large nonlinear programs that can be computationally demanding for conventional CPU-based solvers. Although emerging GPU-based solvers offer substantial computational benefits, their application to process optimization has been limited by the lack of GPU-compatible process modeling tools. We address this gap by prototyping the GPU-compatible process optimization models using an existing GPU-capable optimization software stack, including ExaModels (algebraic modeling system), MadNLP (optimization solver), and cuDSS (linear solver). ExaModels formulates the process optimization problem in a GPU-compatible way by exposing its repeated algebraic structure, while MadNLP and cuDSS solve the resulting nonlinear program on the GPU. This workflow is demonstrated on a $CO_2$ absorber design problem under feed uncertainty, in which a shared column diameter is minimized subject to equilibrium and hydraulic constraints in all scenarios. For the largest case with 5,000 scenarios and 1.5 million variables, the GPU workflow achieves a speedup of approximately 21× over a single-threaded CPU baseline using JuMP, Ipopt, and MA57.

## INTRODUCTION

Large-scale process design optimization is limited by the scalability of current nonlinear optimization workflows. After a flowsheet has been specified, design optimization determines equipment sizes and operating decisions that satisfy process constraints while optimizing economic or performance objectives [1]. When uncertainty in feed conditions, product demand, prices, or operating conditions is present, the design must perform effectively across multiple possible realizations. Scenario-based stochastic formulations address this challenge by assigning scenario-dependent operating variables and process constraints to each uncertainty realization while coupling all scenarios through shared design variables [2, 3, 4]. This improves the robustness of the design, but also causes the size of the problem to grow with the number of uncertainty realizations.

Scalable stochastic process optimization in process systems engineering has advanced through both improved formulations and improved solution methods. On the formulation side, multistage and multiperiod models represent recourse, nonanticipativity, and sequential realization of uncertainty in process design, operation, and planning [3, 4, 5]. On the solution side, decomposition methods, sample-average approximation, and scenario generation and reduction have made these large structured models more computationally tractable [6, 7, 8]. However, fully GPU-resident nonlinear optimization remains unexplored for scenario-based chemical process design problems.

Equation-oriented process models are well suited to GPU acceleration because the constraints exhibit a predictable repeated structure. These models typically consist of material and energy balances, phase-equilibrium relations, summation constraints, operating limits, and process specifications [9]. When repeated across column stages, components, phases, and uncertainty scenarios, these equations generate large sparse nonlinear programs with substantial uniform structure. This structure is well matched to GPU architectures, which are most effective when function and derivative evaluations, as well as linear algebra operations, can be expressed as regular, data-parallel workloads.

Recent advances in GPU nonlinear optimization now make it possible to execute the entire nonlinear programming workflow on accelerator hardware, including model evaluation, automatic differentiation, sparse linear algebra, and nonlinear solver iterations. These capabilities avoid time-consuming host-device data transfers and have produced substantial speedups for large structured nonlinear programs in applications such as optimal power flow [10], but they remain unexplored for scenario-structured chemical process design. In this work, we exploit the repeated structure of equation-oriented process optimization under uncertainty using a fully GPU-resident nonlinear optimization workflow.

The main contributions of this paper are twofold. First, we develop a proof-of-concept GPU-resident implementation of a scenario-based equation-oriented equilibrium absorber design problem using ExaModels.jl and MadNLP.jl [10]. The absorber diameter is a shared design variable, while scenario-dependent operating states satisfy replicated equilibrium process constraints. Second, we demonstrate the scalability of this workflow on problems with up to 5000 scenarios, approximately 1.5 million variables, and 1.05 million equality constraints, achieving a GPU solution time of 33.3 s compared with 695.5 s for JuMP/Ipopt and 707.8 s for the CPU ExaModels.jl/MadNLP.jl configuration. Within the interior-point iterations, the GPU workflow uses CuDSS to solve sparse linear systems, whereas both CPU workflows use MA57. These results demonstrate that GPUs are a viable and scalable platform for scenario-structured process design optimization under uncertainty.

The remainder of this paper proceeds as follows. We first present the deterministic equilibrium process model. We then formulate the scenario-based nonlinear program and describe its ExaModels.jl construction and GPU-resident solution with MadNLP.jl. The case study applies the workflow to an absorber diameter minimization problem under uncertain feed conditions and reports the numerical results. The paper concludes with key findings and directions for future work.

# EQUILIBRIUM PROCESS MODEL

This section defines the equilibrium stage model used for process design under uncertainty. The model consists of Material, Equilibrium, Summation, and Heat (MESH) equations, with auxiliary pressure relations and process specifications.

Let $\mathcal{N} = \{1, \ldots, N\}$ denote the equilibrium stages, numbered from bottom to top, let $\mathcal{C}$ denote the chemical components, and let $\mathcal{R}$ denote the equilibrium reactions.

## Material balances

The material balance conserves component flow across the stage:

$$L_n x_{i,n} + V_n y_{i,n} = L_{n+1} x_{i,n+1} + V_{n-1} y_{i,n-1} + f^{\mathrm{L}}_{i,n} + f^{\mathrm{V}}_{i,n} + \textstyle\sum_{r\in\mathcal{R}} \nu_{i,r}\, \dot{\varepsilon}_{r,n}, \forall i \in \mathcal{C}, n \in \mathcal{N}, \quad (1)$$

where $x_{i,n}$ and $y_{i,n}$ are the liquid and vapor mole fractions of component $i$ on stage $n$, and $x_n = \{x_{i,n}\}_{i\in\mathcal{C}}$ and $y_n = \{y_{i,n}\}_{i\in\mathcal{C}}$ collect the corresponding phase composition vectors; $L_n$ and $V_n$ are the liquid and vapor molar flow rates leaving stage $n$. Thus, $L_{n+1}$ enters from the stage above and $V_{n-1}$ enters from the stage below. Boundary streams, such as $L_{N+1}$ and $V_0$, represent specified feeds or neighboring unit operations. $f^{\mathrm{L}}_{i,n}$ and $f^{\mathrm{V}}_{i,n}$ are optional liquid and vapor component feeds. $\nu_{i,r}$ is the stoichiometric coefficient of component $i$ in reaction $r \in \mathcal{R}$, and $\dot{\varepsilon}_{r,n}$ are the reaction extent rates.

## Equilibrium relations

Phase equilibrium is enforced with equilibrium ratios:

$$y_{i,n} = K_{i,n} x_{i,n}, \qquad \forall i \in \mathcal{C}, n \in \mathcal{N}, \quad (2a)$$

$$K_{i,n} = K_i(T_n, P_n, x_n, y_n), \qquad \forall i \in \mathcal{C}, n \in \mathcal{N}, \quad (2b)$$

where $T_n$ and $P_n$ are the temperature and pressure on stage $n$, $K_{i,n}$ is the equilibrium ratio, and $K_i(\cdot)$ is the thermodynamic correlation that computes the equilibrium ratio from temperature, pressure, and composition.

Reaction equilibrium is written in terms of activities:

$$\textstyle\sum_{i\in\mathcal{C}} \nu_{i,r} \ln a_{i,n} = \ln K^{\mathrm{R}}_{r,n}, \qquad \forall r \in \mathcal{R}, n \in \mathcal{N}, \quad (3a)$$

$$a_{i,n} = a_i(T_n, P_n, x_n, y_n), \qquad \forall i \in \mathcal{C}, n \in \mathcal{N}, \quad (3b)$$

$$K^{\mathrm{R}}_{r,n} = K^{\mathrm{R}}_r(T_n, P_n), \qquad \forall r \in \mathcal{R}, n \in \mathcal{N}, \quad (3c)$$

where $a_{i,n}$ is the activity of component $i$, and $a_i(\cdot)$ is the chosen activity model. Likewise, $K^{\mathrm{R}}_{r,n}$ is the reaction equilibrium constant variable, and $K^{\mathrm{R}}_r(\cdot)$ is its thermodynamic correlation.

## Summation constraints

The summation constraints normalize the liquid and vapor compositions on each stage:

$$\textstyle\sum_{i\in\mathcal{C}} x_{i,n} = 1, \quad \sum_{i\in\mathcal{C}} y_{i,n} = 1, \quad \forall n \in \mathcal{N}. \quad (4)$$

## Heat balances

The heat balance enforces enthalpy conservation. The liquid and vapor molar enthalpies are computed from phase enthalpy models [11]:

$$h^{\mathrm{L}}_n = h^{\mathrm{L}}(T_n, P_n, x_n),\ h^{\mathrm{V}}_n = h^{\mathrm{V}}(T_n, P_n, y_n), \quad \forall n \in \mathcal{N}, \quad (5)$$

where $h_n^{\mathrm{L}}$ and $h_n^{\mathrm{V}}$ are liquid and vapor molar enthalpy, while $h^{\mathrm{L}}(\cdot)$ and $h^{\mathrm{V}}(\cdot)$ are the corresponding enthalpy correlations.

The stage energy balance is

$$L_n h_n^{\mathrm{L}} + V_n h_n^{\mathrm{V}} = L_{n+1} h_{n+1}^{\mathrm{L}} + V_{n-1} h_{n-1}^{\mathrm{V}} + F_n^{\mathrm{L}} h_n^{\mathrm{F,L}} + F_n^{\mathrm{V}} h_n^{\mathrm{F,V}} + Q_n, \quad \forall n \in \mathcal{N}, \quad (6)$$

where $Q_n$ is the heat duty, $F_n^{\mathrm{L}}$ and $F_n^{\mathrm{V}}$ are optional liquid and vapor side feed molar flow rates, and $h_n^{\mathrm{F,L}}$ and $h_n^{\mathrm{F,V}}$ are their molar enthalpies. For boundary stages, $h_{N+1}^{\mathrm{L}}$ and $h_0^{\mathrm{V}}$ denote inlet stream enthalpies.

### Pressure relations

Stage pressures may be specified or computed from pressure drop relations. For adjacent stages $n = 1, \dots, N-1$, one common form is

$$P_n = P_{n+1} + \Delta P_n; \quad P_N = P^{\mathrm{top}}, \quad (7)$$

where $\Delta P_n$ is the pressure drop from stage $n+1$ to stage $n$, and $P^{\mathrm{top}}$ is the specified top stage pressure. Depending on the unit model, $\Delta P_n$ may be specified, set to zero, or computed from a hydraulic correlation.

### Degrees of freedom

For a deterministic problem, (1)–(7) form a square MESH block once the boundary streams, side feeds, heat duty and pressure specifications, as well as any material property specifications are supplied. These specifications close the material, equilibrium, summation, heat, and pressure equations, leaving no degrees of freedom for simulation. In an optimization problem, selected quantities that would otherwise be specified are instead left free and bounded. These quantities become the design and operating degrees of freedom and are optimized subject to the MESH equations.

The deterministic MESH block is sparse because each stage depends only on local variables and adjacent interstage streams.

## SCENARIO-BASED NONLINEAR PROGRAM CONSTRUCTION AND GPU SOLUTION

The MESH block (1)–(7) is embedded in a two-stage, scenario-based nonlinear program (NLP). Let $\Omega$ denote a finite scenario set with probabilities $p_\omega \geq 0$ satisfying $\sum_{\omega \in \Omega} p_\omega = 1$. The deterministic-equivalent stochastic design problem is

$$\min_{\theta, z_\Omega} \quad C^{\mathrm{des}}(\theta) + \textstyle\sum_{\omega \in \Omega} p_\omega \, C_\omega^{\mathrm{op}}(z_\omega; \theta, \xi_\omega) \quad (8a)$$

$$\text{s.t.} \quad F_\omega^{\mathrm{MESH}}(z_\omega; \theta, \xi_\omega) = 0, \quad \forall \omega \in \Omega, \quad (8b)$$

$$G_\omega(z_\omega; \theta, \xi_\omega) \leq 0, \quad \forall \omega \in \Omega, \quad (8c)$$

$$\theta^L \leq \theta \leq \theta^U, \quad (8d)$$

$$z_\omega^L \leq z_\omega \leq z_\omega^U, \quad \forall \omega \in \Omega, \quad (8e)$$

where $\xi_\omega$ collects the uncertain data in scenario $\omega$, including feed flow rates, feed compositions, temperatures, pressures, and other operating parameters. The design vector $\theta$ is chosen before uncertainty is realized, whereas the operating and state vector $z_\omega$ is scenario-specific recourse.

The term $C^{\mathrm{des}}(\theta)$ is the design or capital cost, and $C_\omega^{\mathrm{op}}(z_\omega; \theta, \xi_\omega)$ is the operating cost in scenario $\omega$. The equalities $F_\omega^{\mathrm{MESH}} = 0$ impose the scenario-specific material balances, equilibrium relations, summation constraints, energy balances, and pressure relations defined above. Equation (8c) collects scenario-dependent feasibility requirements, such as hydraulic capacity, pressure, equipment, and process specification limits. The vectors $\theta^L$ and $\theta^U$ are lower and upper bounds on the shared design variables, and $z_\omega^L$ and $z_\omega^U$ are lower and upper bounds on the scenario recourse variables. Relative to the square deterministic MESH block described before, the degrees of freedom in (8) are the elements of $\theta$ and $z_\omega$ that are not fixed as simulation specifications.

### Replicated scenario structure

Equations (8b)–(8c) are replicated for all scenarios. Here, replication means that a common algebraic constraint templaterelat is instantiated across scenarios, stages, components, and reactions, with scenario-specific variables, data, and bounds. The deterministic equivalent is therefore block-separable by scenario, except for the shared design variables, which provide the only cross-scenario coupling. If one deterministic MESH model contains $n_z$ recourse variables and $n_f$ equality constraints, then the scenario model contains $|\Omega| n_z$ recourse variables and $|\Omega| n_f$ MESH equalities, plus the shared design variables and scenario-dependent inequalities. This linear growth with $|\Omega|$ creates the scalability challenge while also providing the exploitable parallel structure.

### Structured model construction with ExaModels.jl

Equation (8) can be implemented in several algebraic modeling systems, including JuMP, AMPL, CasADi, and hand-coded callback interfaces. We use ExaModels.jl, a state-of-the-art platform for sparse NLPs with repeated algebraic structure. Its generator-based syntax declares objective and constraint equations over index sets, exposing the single-instruction, multiple-data structure required for parallel function and automatic differentiation (AD) evaluation on CPUs and GPUs [10].

ExaModels.jl represents each constraint family using a generator: a single algebraic template followed by a `for` clause that specifies the index sets over which the template is instantiated. A representative declaration for the material-balance family is given below.

```
@add_con(core,
    L[n,w]*x[i,n,w] + V[n,w]*y[i,n,w]
  - L[n+1,w]*x[i,n+1,w]
```

```
  - V[n-1,w]*y[i,n-1,w]
  - fL[i,n,w] - fV[i,n,w]
  - sum(nu[i,r]*eps[r,n,w] for r in R)
  for i in C, n in N, w in W;
    lcon = 0.0, ucon = 0.0)
```

In this declaration, the expression defines one scalar material-balance equation, while `for i in C, n in N, w in W` instantiates it over all components, stages, and scenarios. Other MESH constraint families are declared analogously. Increasing the number of components, stages, reactions, or scenarios changes only the index sets and array dimensions, not the algebraic template for each family.

This generator-based formulation exposes the repeated MESH structure directly to the AD and solver interfaces. Manually constructing scalar constraints in an outer Julia loop should be avoided in this formulation because it obscures the common template and limits GPU kernel generation.

### GPU-resident solution with MadNLP.jl

By exposing the parallel sparse NLP structure, the ExaModels.jl implementation above makes the resulting problem suitable for a fully GPU-resident solution with MadNLP.jl. In the GPU configuration, function evaluation, AD, KKT assembly, and sparse linear algebra remain on the device, avoiding repeated host-device transfers during solver callbacks [10].

## CASE STUDY: ABSORBER DIAMETER MINIMIZATION

### Problem setup

Our case study considers the optimal sizing of an absorber column for solvent-based carbon capture, a mature technology that is applied on an industrial scale [12]. Specifically, we aim to find the minimum absorber diameter that ensures feasible, non-flooding operation under all uncertain feed conditions (Fig. 1). For each scenario $\omega \in \Omega$, the inlet gas and liquid flow rates, compositions, temperatures, and absorber operating pressure $P_\omega^{\text{abs}}$ are prescribed. The absorber diameter $d$ is the shared design variable, whereas stage flow rates, compositions, temperatures, and reaction extents are scenario-dependent recourse variables.

The physical absorber is a continuous contacting column. In the optimization model, we approximate it by an equilibrium-stage representation with $\mathcal{N} = \{1, \dots, N_{\text{abs}}\}$ and $N_{\text{abs}} = 10$. The scenario set is generated from seven measured baseline operating cases by sampling random convex combinations, so all sampled feed conditions lie within the convex hull of the measurements. The objective is to minimize the column diameter $d$ that satisfies the gas superficial-velocity limit on every modeled stage

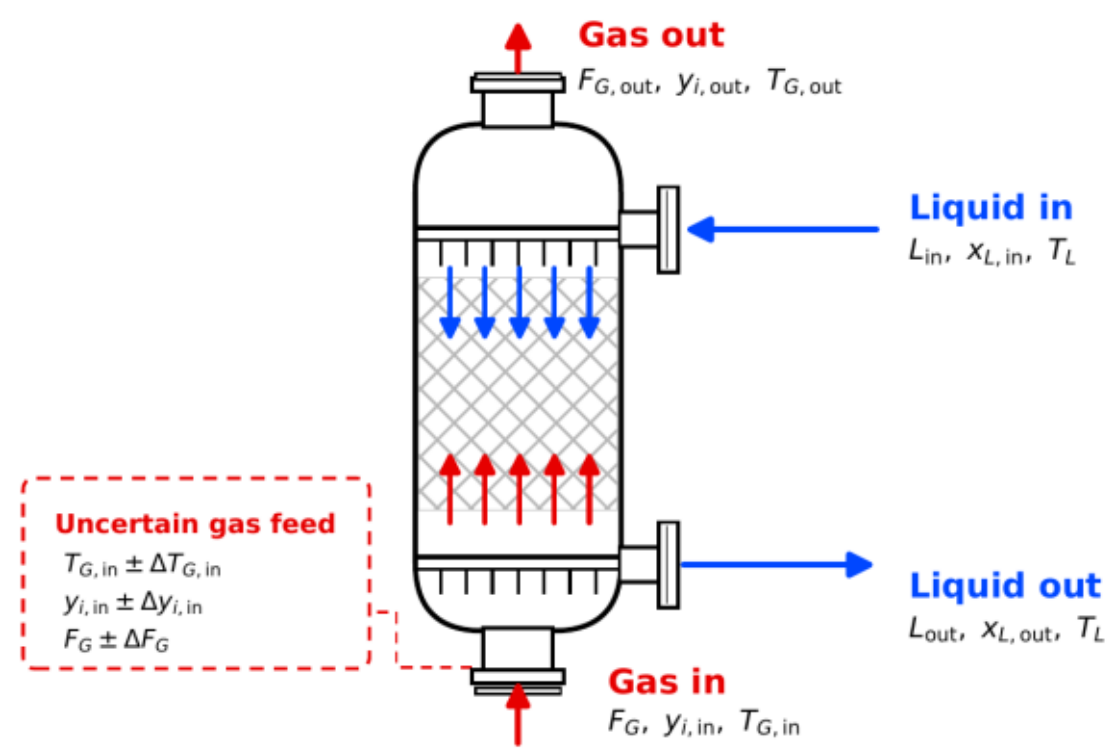

Figure 1. Absorber under feed uncertainty.

in every scenario.

The set of components and the five liquid-phase equilibrium reactions are:

$$\mathcal{C} = \{\text{MEA}, \text{MEAH}^+, \text{H}^+, \text{MEACOO}^-, \text{HCO}_3^-, \text{OH}^-, \text{CO}_3^{2-}, \text{CO}_2, \text{H}_2\text{O}, \text{N}_2, \text{O}_2\}.$$

$$\text{MEAH}^+ \overset{K_1}{\leftrightarrow} \text{MEA} + \text{H}^+ \tag{R1}$$

$$\text{HCO}_3^- \overset{K_2}{\leftrightarrow} \text{H}^+ + \text{CO}_3^{2-} \tag{R2}$$

$$\text{H}_2\text{O} \overset{K_3}{\leftrightarrow} \text{H}^+ + \text{OH}^- \tag{R3}$$

$$\text{MEACOO}^- + \text{H}_2\text{O} \overset{K_4}{\leftrightarrow} \text{MEA} + \text{HCO}_3^- \tag{R4}$$

$$\text{H}_2\text{O} + \text{CO}_2 \overset{K_5}{\leftrightarrow} \text{H}^+ + \text{HCO}_3^- \tag{R5}$$

The temperature-dependent equilibrium constants are computed with the Kent–Eisenberg correlation [13].

$$\ln K_r(T) = A_r + \frac{B_r}{T} + C_r \ln T + D_r T, \quad r = 1, \dots, 5.$$

The vapor phase is modeled as an ideal gas mixture and the liquid phase as an ideal liquid mixture, so vapor fugacities are represented by partial pressures and liquid activities by mole fractions.

Water vapor–liquid equilibrium follows Raoult's law, and $CO_2$ vapor–liquid equilibrium follows Henry's law. For all $n \in \mathcal{N}$ and $\omega \in \Omega$,

$$P_\omega^{\text{abs}} y_{\text{H}_2\text{O},n,\omega} = x_{\text{H}_2\text{O},n,\omega} P_{\text{H}_2\text{O}}^{\text{sat}}(T_{n,\omega}), \tag{9a}$$

$$P_\omega^{\text{abs}} y_{\text{CO}_2,n,\omega} = c_{\text{CO}_2,n,\omega} H_{\text{CO}_2}\left(T_{n,\omega}, x_{\text{MEA},n,\omega}^{\text{app}}, x_{\text{H}_2\text{O},n,\omega}^{\text{app}}\right), \tag{9b}$$

where $P_\omega^{\text{abs}}$ is the prescribed absorber pressure, $P_{\text{H}_2\text{O}}^{\text{sat}}(T)$ is the water saturation-pressure correlation, $c_{\text{CO}_2,n,\omega}$ is the dissolved concentration, and $x_{\text{MEA},n,\omega}^{\text{app}}$ and $x_{\text{H}_2\text{O},n,\omega}^{\text{app}}$ are the apparent liquid compositions [11, 14]. Pressure drop is neglected, so all stages operate at the scenario absorber pressure $P_{n,\omega} = P_\omega^{\text{abs}}$.

The energy balance uses ideal-mixture enthalpies. The liquid and vapor enthalpy flow rates are $\dot{H}^{\mathrm{L}}_{n,\omega} = L_{n,\omega} \sum_{i \in \mathcal{C}_{\mathrm{L}}} x_{i,n,\omega}\, h^{\mathrm{L}}_i(T_{n,\omega})$ and $\dot{H}^{\mathrm{V}}_{n,\omega} = V_{n,\omega} \sum_{i \in \mathcal{C}_{\mathrm{V}}} y_{i,n,\omega}\, h^{\mathrm{V}}_i(T_{n,\omega})$, where $\mathcal{C}_{\mathrm{L}}$ and $\mathcal{C}_{\mathrm{V}}$ denote the liquid and vapor components. The component molar enthalpies are computed from temperature-dependent heat capacities with reference-state corrections [11]:

$$h_i^{\phi}(T) \quad = h_i^{\phi}(T_{\mathrm{ref}}) + \int_{T_{\mathrm{ref}}}^{T} C_{p,i}^{\phi}(\tau)\, d\tau + \Delta h_i^{\phi}(T), \quad \phi \in \{\mathrm{L}, \mathrm{V}\}.$$

The absorber diameter is constrained by the maximum allowable superficial gas velocity. The hydraulic constraint is written as

$$\frac{V_{n,\omega} R T_{n,\omega}}{P_{\omega}^{\mathrm{abs}}} \le u_{\max} \frac{\pi d^2}{4}, \quad n \in \mathcal{N},\ \omega \in \Omega,$$

where $u_{\max}$ is the specified superficial velocity limit and $R$ is the ideal gas constant.

Together with (1), (3a), (4), (6), and (9), the diameter-minimization problem is

$$\min_{d, \{z_\omega\}_{\omega \in \Omega}} \quad d \tag{10a}$$

$$\text{s.t.} \quad \frac{V_{n,\omega} R T_{n,\omega}}{P_{\omega}^{\mathrm{abs}}} - u_{\max} \frac{\pi d^2}{4} \le 0, \quad n \in \mathcal{N},\ \omega \in \Omega, \tag{10b}$$

$$d^L \le d \le d^U, \tag{10c}$$

$$z_\omega^L \le z_\omega \le z_\omega^U, \quad \omega \in \Omega. \tag{10d}$$

## Numerical results

We solved the deterministic equivalent of (10) for $|\Omega| \in \{1, 10, 100, 1000, 5000\}$ scenarios. For the largest instance, with $|\Omega| = 5000$, the resulting nonlinear program contains 1,500,001 variables, 1,050,000 MESH equality constraints, 450,000 variable-fixing constraints, and 50,000 hydraulic inequality constraints. These dimensions arise because each scenario introduces its own recourse variables and an instantiation of the stagewise MESH constraints, while all scenarios share a single first-stage absorber diameter.

All experiments were conducted on a system with dual socket 48-core Intel Xeon Platinum 8568Y+ CPUs (2.3 GHz base, 4.0 GHz max turbo) with 2 TB of system memory and 4x NVIDIA H200 GPUs with 141 GB of VRAM device memory each. For the linear solver, the JuMP+Ipopt and the CPU ExaModels.jl+MadNLP.jl configuration used HSL MA57 v3.11.3; the GPU configuration used cuDSS v0.7.1. CPU runs used one Julia thread and one BLAS thread. All configurations used the same initial point, variable bounds, and solver tolerances of $10^{-8}$. The runs used Julia 1.12.1 and CUDA v6.1.0; the implementations are available at https://github.com/mit-shin-group/FOCAPO_ExaProcess.

We compared three solution setups: CPU-based JuMP/Ipopt, CPU-based ExaModels.jl/MadNLP.jl, and GPU-based ExaModels.jl/MadNLP.jl. The CPU configurations were executed in single-threaded mode to serve as a sequential reference for performance comparison. For small scenario sets, the GPU implementation was not consistently faster because fixed GPU launch and data-management overheads were comparable to the total solution time. As the number of scenarios increased, however, the instantiated MESH equations and hydraulic constraints created enough parallel work across scenarios and stages to amortize these overheads. At $|\Omega| = 5000$, the GPU ExaModels.jl/MadNLP.jl configuration solved the problem in 33.3 s, compared with 695.5 s for JuMP/Ipopt and 707.8 s for the CPU ExaModels.jl/MadNLP.jl configuration (Table 1). These solve times correspond to speedups of $20.9\times$ relative to JuMP/Ipopt and $21.3\times$ relative to the CPU ExaModels.jl/MadNLP.jl configuration. The iteration counts in Table 1 distinguish convergence behavior from per-iteration computational cost. For the large instances, the GPU speedup was primarily attributed to lower per-iteration cost, as the scenario-replicated structure provided sufficient parallel work to exploit the device.

As a sanity check for this prototype problem, the minimum hydraulic margin was within solver tolerance of zero, indicating an active hydraulic constraint and supporting the validity of the optimal solution.

# CONCLUSION

This work presented a proof-of-concept, fully GPU-resident workflow for scenario-based equation-oriented process design under uncertainty. ExaModels exposes replicated MESH constraint structure, while MadNLP and cuDSS solve the resulting nonlinear program on the GPU. The absorber prototype problem demonstrated that GPU acceleration becomes beneficial once the scenario set is large enough to amortize fixed overheads, achieving an order-of-magnitude speedup relative to single-threaded CPU configuration at the largest tested problem. Future work will extend the workflow to full-detail absorber and flowsheet models and further exploit scenario block structure in GPU-based linear algebra and solution of the KKT systems [15].

# ACKNOWLEDGEMENTS

DS, BH, and SS received support from the Future Energy Systems Center through the MIT Energy Initiative. This publication is based on work supported by the U.S. Department of Energy under contract DE-AC02-06CH11357. Funding was provided by the U.S. DOE Hydrocarbons and Geothermal Energy Office under FWP-DOE-145-SE-14205HPC4E-25.

# DECLARATION OF USE OF AI

Generative AI tools were used solely to assist with language editing and to improve the clarity, concision,

Table 1. Solve times and interior-point iteration counts for each solution setup.

| $\lvert\Omega\rvert$ | JuMP/Ipopt/MA57 | | ExaModels/MadNLP/MA57 (CPU) | | ExaModels/MadNLP/CuDSS (GPU) | |
|---|---|---|---|---|---|---|
| | Solve (s) | Iterations | Solve (s) | Iterations | Solve (s) | Iterations |
| 1 | 0.03 | 18 | 0.008 | 20 | 0.22 | 16 |
| 10 | 0.32 | 29 | 0.11 | 26 | 0.62 | 27 |
| 100 | max iter exceeded | | 2.61 | 46 | 0.83 | 42 |
| 1000 | 91.9 | 52 | 369.4 | 599 | 15.8 | 78 |
| 5000 | 695.5 | 80 | 707.8 | 115 | 33.3 | 139 |

and organization of the manuscript. All suggestions were reviewed and edited by the authors, who take full responsibility for the content. AI was not used to produce research results, models, data, or numerical experiments. AI is not an author of this work.

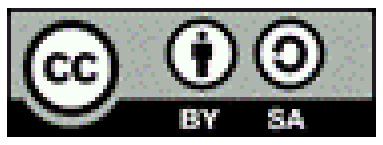